\documentclass[10pt]{article}

\usepackage[centertags]{amsmath}
\usepackage{amsfonts,amssymb,amsthm}
\usepackage{authblk}
\usepackage{caption,graphicx}
\usepackage[text={6.0in,8.6in},centering,letterpaper]{geometry}
\usepackage{hyperref}
\usepackage{mathrsfs}
\usepackage{mathtools}
\usepackage[numbers,comma,square,sort&compress]{natbib}
\usepackage{subcaption}
\usepackage{tikz}
\usepackage{verbatim}
\usepackage{url}
\newtheorem{remark}{Remark}

\DeclarePairedDelimiter{\nrm}\lVert\rVert
\DeclarePairedDelimiter{\ip}\langle\rangle

\newcommand{\R}{\ensuremath{\mathbb{R}}}
\newcommand{\e}{\ensuremath{\varepsilon}}
\newcommand{\p}{\ensuremath{\partial}}

\newcommand{\s}{\ensuremath{\sigma}}

\numberwithin{equation}{section}
\newcommand{\sref}[1]{(\ref{#1})}                       

\begin{document}
\author{C.H.S. Hamster}
\affil{Dutch Institute for Emergent Phenomena,\\University of Amsterdam.\\
Email:  {\normalfont{\texttt{c.h.s.hamster@uva.nl}}}}
\title{Stochastic Nonlinear Waves in all the Right Phases}
\date{\today}
\maketitle

\begin{abstract}
    I present an overview of different phases used in the literature to describe the phase of travelling waves in reaction-diffusion equations. The discussion is tailored to understand the effectiveness of these different phase descriptions in a stochastic setting in 1D. The explicit computation of the autonomous isochronal phase for the Nagumo equation is a new result. 
\end{abstract}

\section{Introduction}
In the recent decade, efforts to study nonlinear waves in a stochastic setting have picked up pace. As the field is relatively new, there are no established practices, which is explicitly shown in a still-developing zoo of different phase descriptions for stochastic travelling waves. Of course, this problem is not induced by the stochasticity; already in the classic deterministic literature, many different phases can be found. In the context of Burgers' equation, Zumbrun showed that a large family of phases are equivalent~\cite{zumbrun2009}. This feels intuitive, as we expect stability of travelling waves to depend on, e.g., the chosen normed vector space, but not on the specific choice of location. 

The goal of this work is to present an overview of the different approaches taken in the literature and to discuss the similarities and differences, and the strengths and weaknesses of each method, specifically in the stochastic setting. Conversations with colleagues convinced me that such an overview is necessary, as the differences between the different descriptions can be subtle, naming conventions are not very strict, and it is not always possible to connect a phase description to a specific author. Two remarks: First, this list is non-exhaustive and comments are always welcome. Secondly, the different phases are ordered thematically, and within each theme, they are ordered randomly. It is therefore purely incidental that the first phase described is my own. 

Before we study all the different phases, maybe we should pause and discuss why phase descriptions are essential in the first place, both deterministically and stochastically. All of the equations discussed below are translation invariant, and therefore, there is no unique travelling wave, but a 1D family of travelling waves indexed by the position. In the stability analysis of the travelling wave, the translation invariance always results in a 0 eigenvalue for the linearisation around the travelling wave, corresponding to the eigenvector pointing along the family of travelling waves. The fact that the dynamics along the family of travelling waves is neutral, i.e. not decaying exponentially back to the wave or growing exponentially away from it, means that the noise has a large influence in this direction. When the noise is translation invariant, the translation invariance remains a symmetry of the SPDE. When an equation has more than one symmetry, such as the Korteweg-de Vries equation, which has two, this needs to be reflected in the phase~\cite{westdorp2024long}.

The sensitivity of travelling waves with respect to the neutral direction is already present in the deterministic setting. Any perturbation at $t=0$ that is not orthogonal to the neutral direction results in a shift of the wave. Therefore, it does not make sense to talk about the stability of a single travelling wave, but we talk about the stability of the family of travelling waves. The fact that a perturbation causes the solution to converge to a shifted version of the wave is called \textit{orbital stability}. In a stochastic setting, orbital stability is often not possible as the wave is not a stationary solution of the SPDE, but we do observe numerically that the stochastic wave remains close to the deterministic manifold of travelling waves. Hence, we could refer to this phenomenon as stochastic meta-orbital stability, but this is not a term that is used in the literature.

To set the stage, we will consider mostly equations of the following form:
\begin{align}
    du=[Du_{xx}+f(u)]dt+\sigma g(u)dW^Q_t, \hspace{3mm} u(t)\in\mathbb{R}^n \hspace{3mm} x\in\mathbb{R}, \hspace{3mm} t\in\mathbb{R}^+.
\end{align}
That is, for now, we assume one dimension in space, but $u$ can have multiple components. For $\sigma=0$, we assume that the equation has a well-studied travelling wave. Famous examples are the Nagumo equation\footnote{Also named Allen-Cahn, Huxley, Slogl, Fisher with a strong Allee effect or just bistable RDE.} and the FitzHugh-Nagumo equation.  

Next, let us discuss $W^Q_t$. This can be just Brownian motion, but also an infinite-dimensional process. Often, it is assumed that $W^Q_t$ is a $Q$-Wiener process, and explicitly not a cylindrical $Q$-Wiener process, meaning that $Q$ is of trace class. On $\R$, this implies that the kernel $q$ of $Q$ is not translation invariant. To what extent this is a problem will be discussed in the next sections, but this is explained in greater detail in the introduction of my thesis~\cite{HamsterThesis}. 

\paragraph{Notation}
Of course, all authors mentioned here use their own notation. I can't streamline everything, so for clarity, I will stick to the following notation. Deterministic solutions are written in lower case, stochastic solutions in capitals. I denote the deterministic travelling wave with $(\Phi_0,c_0)$, to differentiate it from perturbed versions.  

The variable $\xi$ is always used to indicate the travelling wave coordinate $\xi=x-c_0t$. This choice implies that the linearization around the travelling wave has the following form:
\begin{align}
    \mathcal{L}v=Dv''+c_0v'+f'(\Phi_0)v.
\end{align}
Other authors, e.g. Sattinger~\cite{Sattinger}, use a different sign convention. In general, any shift gets a minus sign, i.e. we say that a wave in the coordinate $\xi$ is shifted by $c_0t$, not by $-c_0t$. In \cite{Inglis}, this is the other way around.

As is well known, the linear operator above has a 0-eigenfunction $\Phi_0'$. Furthermore, the adjoint operator $\mathcal{L}^*$, defined as 
\begin{align}
    \mathcal{L}^*v=Dv''-c_0v'+f'(\Phi_0)v,
\end{align}
also has a zero eigenfunction $\sim e^{c_0\xi}\Phi'_0$. By $\psi$, I denote the normalized adjoint eigenfunction, i.e. $\ip{\Phi'_0,\psi}_{L^2(\R)}=1$. 

\paragraph{What to shift?}
Note that there are two ways to introduce a new coordinate system. For solution $u(t)$ and phase $\gamma(t)$, we could either choose $\gamma$ such that $\nrm{u(\cdot,t)-\Phi(\cdot\pm\gamma(t))}\to0$ or \\$\nrm{u(\cdot\pm\gamma(t),t)-\Phi(\cdot)}\to0$ (the $\pm$ depends on your preference). Of course, these two options are completely equivalent, but when it comes down to proving things, there is a major difference. The first option will lead to a linearization around the travelling wave $\Phi(\cdot+\gamma(t))$, while the second option leads to a linearization around $\Phi(\cdot)$. Now, for the spectral properties of the linear operator, this difference is moot, but when closing a nonlinear iteration argument to prove stability, the time-dependence of $\Phi(\cdot+\gamma(t))$ (and hence $\mathcal{L}$) becomes important\footnote{In~\cite{kapitula}, this $\gamma$-dependence is referred to as `irksome'. }, see for example~\cite{Inglis}. However, the advantage is that all derivatives coming from the (stochastic) chain rule are on the (smooth) shifted wave. When the solution is shifted, all derivatives are on the shifted solution, which need to be controlled in a proof. Especially in a stochastic context, controlling derivatives can be tricky.    

As an interesting remark, the deterministic versions of several approaches mentioned here are not readily available in the literature. The phase-tracking approach by Hamster and Hupkes was pioneered by Howard and Zumbrun to study shockwaves~\cite{zumbrunhoward}, a much harder problem than the stability of travelling waves in the Nagumo equation. See also the lecture notes by Margaret Beck\footnote{\url{https://math.bu.edu/people/mabeck/bremen-lecture-notes-final-2016.pdf}}. As far as I know, nobody went back to apply it to the easier case of the bistable reaction-diffusion equation. For the case of Inglis and MacLaurin~\cite{Inglis}, the deterministic proof does not exist (I assume) because the proof does not work for travelling waves, as we will discuss below. The formulation of the predicted phase by Van Winden in \cite{vanWinden2024noncommutative} does resemble classic nonlinear iteration, but the general formulation in Banach spaces is new even for the deterministic case. Stannat studied both the deterministic and stochastic versions of this approach, as the deterministic version was not described in the literature~\cite{stannatnag}. 

\paragraph{Naming convention}
I always referred to my methods as `phase tracking', or slightly more extended, `phase tracking \`a la Zumbrun'. However, the term `stochastic freezing', after the deterministic approach initiated by Beyn~\cite{beyn2004freezing}, is also apt. However, these terms only apply to the fact that I shift the solution, not the wave. A more complete name would be something like ``stochastic freezing on the adjoint eigenspace." For the approach developed by Stannat et al.~\cite{lang,stannatkruger}, the phase lags behind the `true' phase, hence the term phase-lag approach. Cartwright and Gottwald~\cite{cartwright2019,cartwright2021collective} use the term `collective coordinates', borrowed from a long history in physics. The approach by Mueller et al.~\cite{mueller1995,mueller2011} depends on level sets; hence, just referring to this phase as a (stochastic) level set suffices. Adams and MacLaurin provided us with the term isochronal phase, and this term/method itself has a long history in the study of oscillators~\cite{winfree1974patterns,guckenheimer1975isochrons}. 

The trickier name is the variational phase, as used in~\cite{Inglis}. The term `variational' comes from the fact that a phase is chosen that minimises the $L^2$-difference between the shifted travelling wave and the solution, i.e. it solves a variational problem. However, as we will discover, the combination of the variational phase with the proof techniques as applied in \cite{Inglis} only works for problems where the linearization around the wave is self-adjoint. When we formulate the problem in a weighted space that makes the linearization self-adjoint, the variational phase does work, as was done in MacLaurin and Bressloff~\cite{maclaurin2020}. The issues with the proof in \cite{Inglis} were solved by MacLaurin in~\cite{maclaurin2023phase} by also projecting onto the adjoint eigenspace, but he still calls the phase variational, even when there is no variational problem underlying the new phase. Can we give a name to this approach? In fact, the method described there is a stochastic version of the deterministic proof in~\cite{kapitula}, which again goes without a name.\footnote{Promislow prefers the term `modulational stability'.} I believe that the problem here is that not-freezing is the default, and therefore no standard is available, apart from `co-moving'.  Hence, we can still, as a community, converge on an appropriate name (and it could very well be possible that that name would be `variational phase'). A quite verbose option would be ``stochastic co-moving with the adjoint eigenspace".   

\paragraph{Not in these notes} 
In this work, I will focus mainly on 1D equations. However, theory on stochastic patterns in 2D has also been developing, see for example \cite{bosch2024multidimensional,ZucalRotatingWaves,kamphuis2025microscopic,chiusole2025existence}. Furthermore, I'll mainly focus on parabolic reaction-diffusion equations. Hence, dispersive equations like Korteweg-de Vries are mentioned, not studied, even though results are available \cite{westdorp2024long,cartwright2021collective}. Many authors mentioned here also worked on neural-field equations, which again will be mentioned but not studied in detail. For results on stochastic neural-field equations, see e.g. \cite{bressloff,stannatkrugerNF,lang,agathe2025stability}. In this work, I always take the viewpoint that I want to understand the average dynamics of the wave, i.e. the average taken over many different realizations of the noise in the SPDE. This is in contrast to the random dynamical system approach, where the SPDE is studied for one realization of the noise, see for example \cite{kuehn2025synchronization}. Parallel to the papers discussed in this overview, the study of sharp interfaces in stochastic versions of equations like Allen-Cahn or Cahn-Hilliard received considerable attention \cite{funaki1995scaling,antonopoulou2016motion,weber2014sharp}. In this approach, both the steepness of the interface and the intensity of the noise depend on the same small parameter and the dynamics is subsequently studied in the limit of this parameter to 0. The connection to the other approaches here should still be explored. Also not mentioned here is the onslaught of scam papers on stochastic travelling waves. After all, the name of the document is ``stochastic nonlinear waves in all the right phases", not ``wrong phases". For an overview of wrong phases, see \cite{hamster2025right}. 

\paragraph{History in physics}
One of the oldest papers in the math literature on stochastic travelling waves is \cite{mueller1995}. However, it would take another 15 years for the subject to gain steam with papers like \cite{mueller2011,lord2012} and \cite{bressloff}. However, this subject has several decades of history in the physics literature.
By the nature of physics, these computations are formal and ad hoc. This can have the advantage of flexibility, for example in~\cite{Armero1996}, where the Nagumo equation is studied in all 3 regimes (bistable and monostable (linear and nonlinear)) at the same time, while all results presented here, apart from Mueller et al.~\cite{mueller1995,maclaurin2023phase}, are for the bistable case (even though most could be extended to the monostable-nonlinear case by moving into weighted spaces). 

The results in~\cite{Armero1996} nicely capture the overall behaviour of the wavespeed as a function of the noise intensity, but deviate from the results in~\cite{Hamster2020} when zoomed into finer detail. Luckily, the authors already warned us that their expansion is incomplete. Furthermore, note that the deviations in~\cite{Armero1996}, and the physics literature in general, purely come from the It\^o-Stratonovich correction, as we will show. 

Suppose we ask the following question: Given a phase description $\Gamma_\sigma(t)$, what is the lowest order expansion in $\s$? In~\cite{garciaspatiallyextended} was shown that at the lowest order in $\sigma$, the phase can be described by a scaled Brownian motion, and the variance is given by
\begin{align}
\label{eq:varphase}
    \ip{g(\Phi_0)Qg^T(\Phi_0)\psi,\psi}t.
\end{align}
An earlier version was already described in~\cite{pasquale1992} in 1992. 
In my opinion, any practical phase description should be able to recover this expression. However, this expression is only valid for the physically relevant case of translation-invariant noise. Most of the works cited here will assume that $Q$ is of trace class and hence not translation-invariant and therefore the expression above is not recovered. This is not an issue with the definition of the phase, but purely with the mathematical techniques used. 
Note that most of the technical discussions in this document on subjects like ``Is $Q$ of trace class or not" or ``Which weight do we need to make a linear operator self-adjoint in a weighted space" become moot when you never formulate your problem in terms of Hilbert spaces in the first place and stick to integrals.

\paragraph{Organization}
This overview is split into three major parts. First, a part where the solution is shifted instead of the wave, then a part where the waves are shifted instead of the solution. Next, a short section with some miscellaneous approaches. The third major part is a study of the Isochronal phase and the explicit computation of the autonomous isochronal phase for the Nagumo equation, which shows that it agrees with the second-order speed corrections from the Hamster and Hupkes approach. 

\section{Shift the solution}
In this section, we study phases where the solution is shifted instead of the wave. In the literature, this approach is sometimes known as \textit{freezing}, after the works by Beyn and Th\"ummler~\cite{beyn2004freezing}, or as phase tracking, after the pioneering work in understanding shockwaves by Howard and Zumbrun~\cite{zumbrunhoward}. I would like to stress here again that, in principle, the choice of shifting the wave or the solution is independent of the choice of phase. However, as the choice of shift influences the proof techniques, I decided to keep this splitting. 

\subsection{Hamster and Hupkes: Adjoint eigenfunction}
\label{sec:HH}
I unapologetically start with my own approach. For educational purposes, let us start with the deterministic equation with one component:
\begin{align}
    u_t=Du_{xx}+f(u), \hspace{3mm} u(t)\in\mathbb{R} \hspace{3mm} x\in\mathbb{R}, \hspace{3mm} t\in\mathbb{R}^+.
\end{align}
We introduce a phase
\begin{align}
    \gamma(t)=\gamma_0+ct+\int_0^ta(u(s))ds
\end{align}
together with the splitting
\begin{align}
    u(\cdot+\gamma(t),t)=\Phi+v(t).
\end{align}
This results in the following equation for $v(t)$:
\begin{align}
    v_t=D(v_{xx}+\Phi'')+(c+a(\Phi+v)(v_x+\Phi')+f(\Phi+v).
\end{align}
Hence, $v$ is what remains when we shift the full solution $u$ with $\gamma$ and subtract $\Phi$. At this point, we are still agnostic about whether or not $(\Phi,c)=(\Phi_0,c_0)$.
We will choose $(\Phi,c)$ and $\gamma(t)$ in order to satisfy one (or, as it turns out, both) of the following goals. \\
\textbf{Goal 1:} \textit{make $\gamma(t)-ct$ as small as possible by choosing $(\Phi,c)$ and $a$.} Ideally, we remove all constant and linear terms from $a(v)$ in $v$. Hence, let us now assume that $a(0)=\partial_va(v)=0$, which we have to check later. 

As we have introduced a one-dimensional degree of freedom, we can implement a one-dimensional constraint. As we work in Hilbert spaces, we assume there is a function $h$ such that $\ip{v(t),h}=0$. This implies
\begin{align}
0&=\frac{\p }{\p t}\ip{v,h}_{L^2(\R)}.
\end{align}
Upon linearising around $\Phi$, we then get
\begin{align}
\label{eq:defav}
a(v)&=-\frac{\ip{D\Phi''+c\Phi'+f(\Phi),h}_{L^2(\R)}+\ip{\mathcal{L}v,h}_{L^2(\R)}+\ip{N(v),h}_{L^2(\R)}}{\ip{\Phi'+v',h}_{L^2(\R)}}.
\end{align}
Hence, if we indeed want a nonlinear $a$, we first must remove the constant term, i.e. we must choose $(\Phi,c)=(\Phi_0,c_0)$. Secondly, we can remove the linear term by ensuring $\ip{\mathcal{L}v,h}_{L^2(\R)}=0$, which implies that $h$ must be in the kernel of the adjoint linear operator $\mathcal{L}^*$, i.e. we choose $h=\psi$. Note that $a(v)$ is now indeed nonlinear in $v$, as we assumed. Doing the computation without immediately choosing $(\Phi,c)=(\Phi_0,c_0)$ is quite convoluted here, but this way of thinking turns out to be useful for the stochastic case. 

\textbf{Goal 2:} \textit{set up a framework that allows us to use semigroup theory.} When we, for convenience, immediately start at $(\Phi_0,c_0)$, we get
\begin{align}
    v_t=\mathcal{L}v+a(v)\partial_x(v+\Phi_0)+N(v).
\end{align}
In mild form, this becomes
\begin{align}
\label{eq:MildSol}
    v(t)=S(t)v_0+\int_0^tS(t-s)[N(v(s))+a(v(s))\partial_x(\Phi_0+v(s))]ds,
\end{align}
where $S(t)$ is the analytic semigroup generated by $\mathcal{L}$. 
Under the assumption that we have a sectorial linear operator with a spectral gap and an isolated simple eigenvalue at 0, the semigroup applied to the nonlinear term decays exponentially, as long as the projection onto the neutral mode $\Phi_0'$ of the semigroup is 0. Here it is important to note that the projection of a function $v$ onto the neutral mode of the semigroup is not $\ip{v,\Phi_0}_{L^2}\Phi_0$, but $\ip{v,\psi}_{L^2}\Phi_0$, see Ch. 4 in~\cite{kapitula}. Upon choosing the integrant of the integral above orthogonal to $\psi$, we arrive again at Eq.~\sref{eq:defav}. Hence, both goals result in the same phase descriptions. 

When we now wish to prove orbital stability, we only have to prove that solutions to Eq.~\sref{eq:MildSol} decay exponentially to zero, which in turn implies that $\gamma(t)-c_0t$ converges to a finite value $\gamma_\infty$. This can be done using a standard nonlinear iteration argument. Note that this result is independent of the choice of phase, we could have chosen any $L^2$-function as a reference function, or even $\delta$-functions if you prefer studying level sets. However, the proof with semigroups only works when we choose $\psi$. We did get orbital stability for free for the description with the variational phase as the norm of $v_\mathrm{var}(t)$ in the variational phase is, by definition, smaller than the norm of $v(t)$.

\subsubsection{Generalization to stochastic case}
We started in the deterministic case with an unknown $(\Phi,c)$ and an unknown function $h$ that defined the phase. Can we do the same in the stochastic case? We define a stochastic phase $\Gamma(t)$ such that 
\begin{align}
    V(t)=U(\cdot+\Gamma(t),t)-\Phi_\s
\end{align}
remains orthogonal to a function $h_\sigma$, where $\Gamma$ is defined as
\begin{align}
    \Gamma_\s(t)=\Gamma_0+c_\sigma+\int_0^ta_\s(V(s))ds+\sigma\int_0^tb_\s(V(s))dW_t^Q.
\end{align}
In an ideal world, we would follow these steps, a generalization of Goal 1 above:
\begin{enumerate}
    \item Choose some reference function $h_\s$ (in $L^2$).
    \item Compute $b_\s$ such that the stochastic part of the equation for $V$ remains orthogonal to $h_\s$.
    \item Compute $(\Phi_\s,c_\s)$ such that the deterministic part disappears when $V=0$, i.e. $a_\sigma(0)=0$. Note that $(\Phi_\s,c_\s)$ depends on $h_\s$. 
    \item Linearise the equation around $\Phi_\s$. 
    \item Choose $h_\s$ as the adjoint eigenfunction of the linear operator computed in step 4, i.e. $\partial_Va_\sigma(0)=0$.  
\end{enumerate}
However, in practice, these steps are impossible, as all the computations depend on $h_\s$, which is only determined in step 5. Furthermore, computing the linearization involves computing the Fr\'echet derivatives of $a_\s(v)$ and $b_\s(v)$; hence, the existence of a 0-eigenvalue cannot be deduced by applying $\partial_\xi$ to the equation for $(\Phi_\s,c_\s)$. Hence, we opted for simplicity to choose $h_\s=\psi$ and only follow steps 1-3. With this choice, the phase $\Gamma_\s(t)$ does not become fully nonlinear, but the linear parts are of $\mathcal{O}(\s^2)$. However, it is important that we defined the phase $\Gamma_\s(t)$ with respect to the reference function $(\Phi_\s,c_\s)$, not with respect to $(\Phi_0,c_0)$!
\begin{remark}
Some explanation of the word `important' in the last sentence is necessary. It is important to realise that this choice is made in the four Hamster and Hupkes papers. I also believe that this choice is convenient. However, it is in no way necessary. We could just have defined a phase $\Gamma_0(t)$ with respect to $(\Phi_0,c_0)$. This would imply that $a(0)\neq0$, but note that $a(0)=\mathcal{O}(\s^2)$. Hence, in the Hamster and Hupkes approach, the average wavespeed at $\mathcal{O}(\s^2)$ gets contributions from $c_\s$ and $\lim_{t\to\infty}t^{-1}E[\Gamma_\s(t)]$, while in the other case, we would only have to study  $\lim_{t\to\infty}t^{-1}E[\Gamma_0(t)]$. If done correctly, both approaches should result in the same expansions. In \cite{BoschHamsterHupkes}, $(\Phi_0,c_0)$ is chosen as the reference function.  
\end{remark}

\subsubsection{Pros and Cons}
This phase description has a very clear con in the stochastic case: the equations become horrible. In terms of the number of complicated terms, there is not much difference with the not-frozen\footnote{Melted?} approach. However, the main problem is that all derivatives are now on the solution $V$. Hence, the stability proof needs to deal with the $V_{xx}$ and $V_{x}$ terms, which luckily can be done, but this is far from trivial. The main advantage of this approach is that it leads to straightforward Taylor expansions in $\sigma$ that can approximate the solution to any order in $\sigma$ (up to a stopping time). We can define the limiting objects
\begin{align}
    \begin{split}
        c_\mathrm{lim}&=\lim_{t\to\infty}E[t^{-1}\Gamma_\s(t)],\\
        \Phi_\mathrm{lim}&=\lim_{t\to\infty}E[U(\cdot+\Gamma_\s(t))].
    \end{split}
\end{align}
Now, these two objects cannot exist, as the solutions are only defined up to a stopping time, but for each term in the Taylor expansion for $\Gamma$ and $V$, these limits do exist, and overlap with numerically computed averages. These expansions can be made rigorous~\cite{BoschHamsterHupkes}. Another way to think about $(\Phi_\mathrm{lim},c_\mathrm{lim})$ is in terms of quasi-stationary distributions. We will explore this more in the section on the isochronal phase. An intuitive explanation can be given for why the approximations of the object $(\Phi_\mathrm{lim},c_\mathrm{lim})$ overlap nicely with the numerical simulations. The convergence towards the quasi-stationary distribution appears to be exponentially fast, while we have to wait exponentially long for the solutions to hit the stopping times.

\subsection{Lord and Th\"ummler: numerical freezing on a reference function.}
In~\cite{lord2012}, the authors developed numerical techniques to study the properties of stochastic waves in a stochastic reaction diffusion equation of the form
\begin{align}
    du=[u_{xx}+f(u)]dt+g(u)dW^Q_t,
\end{align}
both in It\^o and Stratonovich interpretation. 
To reduce the domain size needed for computations, the solution is frozen using a phase condition. The chosen phase condition is that the deviation of the frozen solution from a reference function $\hat u$ should be orthogonal to $\partial_x\hat u$, as this minimizes the deviations in the $L^2$ norm. The choice of what constitutes a `good' reference function now becomes a purely numerical question. For example, taking $\hat u$ to be a Heaviside function leads to numerical instability. When we take the most straightforward option $\hat u=\Phi_0$, this phase is identical to the approach taken in \cite{BoschHamsterHupkes}. 

\section{Shifting the wave}

\subsection{Inglis and MacLaurin: variational approach.}
In this article~\cite{Inglis}, both reaction-diffusion equations and neural field equations are studied. The variational phase is defined as follows. For each time $t$, $\gamma(t)$ is defined\footnote{Warning: Inglis and MacLaurin define $\Phi_\gamma$ as $\Phi(\cdot+\gamma)$, while I define $T_\gamma\Phi=\Phi(\cdot-\gamma)$.} as the phase that minimizes the difference between $u(t)$ and $\Phi_0(\cdot +\gamma(t))$, i.e. we have the following variational problem: 
\begin{align}
    \gamma(t)=\inf_{\alpha\in \R}\nrm{u(t)-\Phi_0(\cdot+\gamma(t))}_{L^2}^2.
\end{align}
If this global minimiser exists, it can be found by computing
\begin{align}
    \frac{d}{d\gamma(t)}\nrm{u(t)-\Phi_0(\cdot +\gamma(t))}^2=-2\ip{u(t)-\Phi_0(\cdot +\gamma(t)),\Phi'_0(\cdot +\gamma(t))}=0.
\end{align}
Hence, the variational approach implies that we must choose the phase such that deviations from the travelling wave are orthogonal to the shifted eigenfunction. This approach can be directly extended to a stochastic setting. The authors derive an SDE for the variational phase, which in turn allows them to write down an SPDE for the deviation $Z(t)=U(t)-\Phi_0(\cdot+\Gamma(t))$. The goal is now to show that $Z(t)$ remains small (in some proper sense) using semigroup theory. However, in Goal 2 in section \ref{sec:HH}, I stated that the adjoint eigenfunction is necessary to use semigroup theory, while~\cite{Inglis} also uses semigroup theory. Is there a contradiction here? No, as the variational approach only works for self-adjoint linear operators. The authors cite the classic Volpert, Volpert and Volpert~\cite{volpert1994traveling} to conclude that, under the same spectral assumptions as we used, a semigroup $S(t)$ can be split as
\begin{align}
    S(t)v=V(t)v+Pv
\end{align}
where $P$ is \textit{the} projection onto $\Phi_0'$ and $V(t)$ satisfies the bound
\begin{align}
    \nrm{V(t)v}\leq e^{-bt}\nrm{v},
\end{align} 
for some $b>0$. 
There are two problems here. Firstly, the bound above, as found in~\cite{volpert1994traveling} is
\begin{align}
    \nrm{V(t)v}\leq Me^{-bt}\nrm{v}.
\end{align} 
In~\cite{Inglis}, the assumption that $M=1$ is essential to close the proof, but this assumption is not explicitly made. Secondly, and this is more fundamental, in~\cite{Inglis}, there is no explicit mention of what $P$ is, while in~\cite{volpert1994traveling}, it is explicitly stated that
\begin{align}
    P=\frac{1}{2\pi i}\oint_{C_0}R(\lambda;\mathcal{L})d\lambda ,
\end{align}
i.e. the contour integral of the resolvent around the pole at 0, 
which simplifies in this setting to
\begin{align}
    Pv=\ip{v,\psi}\Phi_0'.
\end{align}
The implicit assumption made in~\cite{Inglis}, is that $P$ is given by
\begin{align}
    Pv=\ip{v,\Phi_0'}\Phi_0'.
\end{align}
As $P$ is unique, this cannot be correct. Note that $\Phi_0'=\psi$ only when $\mathcal{L}$ is self-adjoint, i.e. for standing waves when $c_0=0$. Hence, the exponential bounds in~\cite{Inglis} on the semigroup cannot be deduced from~\cite{volpert1994traveling}. Given the fact that all phase descriptions of this type are equivalent in the deterministic case, see ~\cite[App. C]{zumbrun2009}, we know that the solution will still decay exponentially, but this cannot be shown directly using semigroups. 

\begin{figure}
\centering
\begin{subfigure}{.49\textwidth}
  \centering
  \includegraphics[width=1\columnwidth]{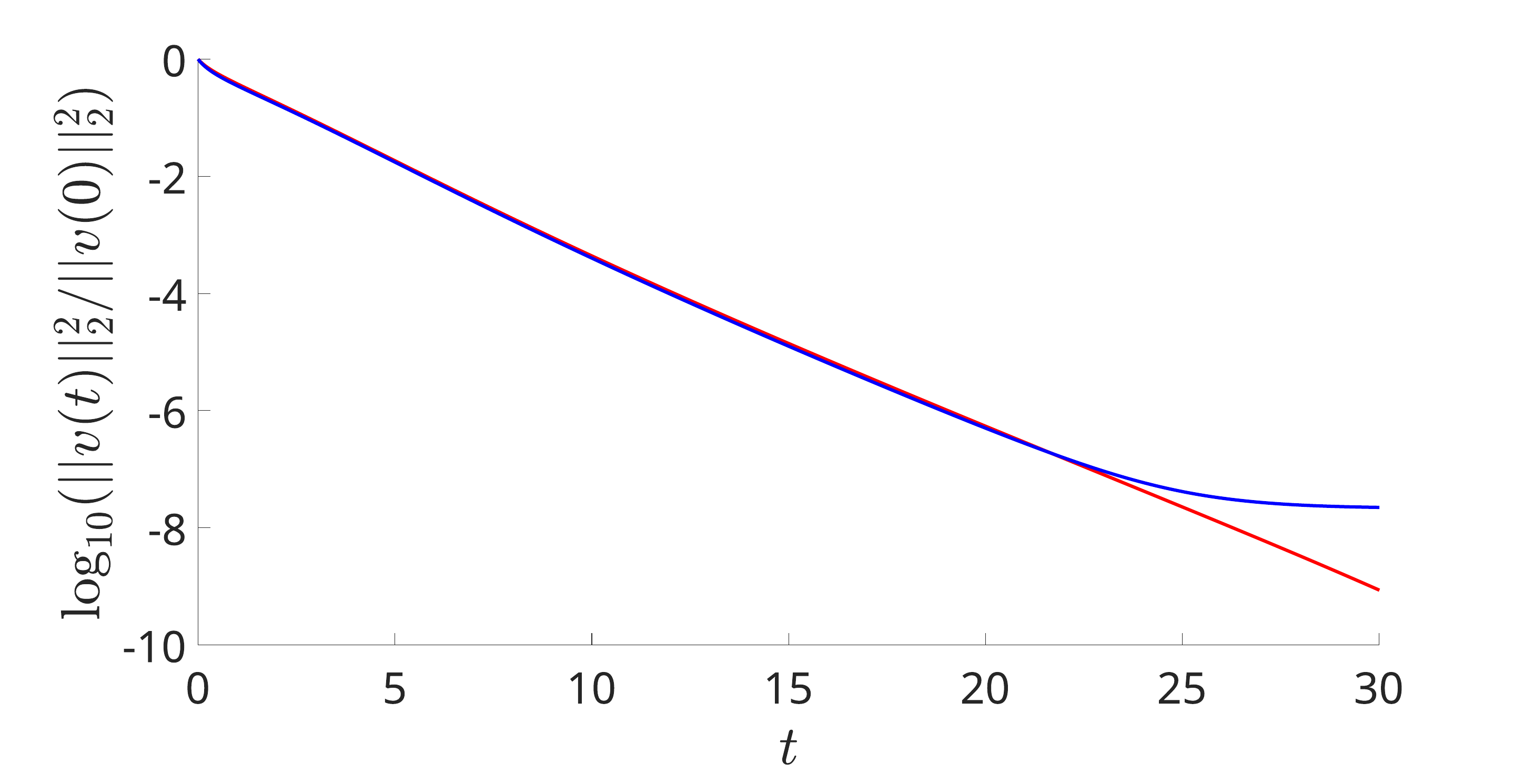}
  \caption{}
      \label{fig:DecayA}
\end{subfigure}
\begin{subfigure}{.49\textwidth}
  \centering
\includegraphics[width=1\columnwidth]{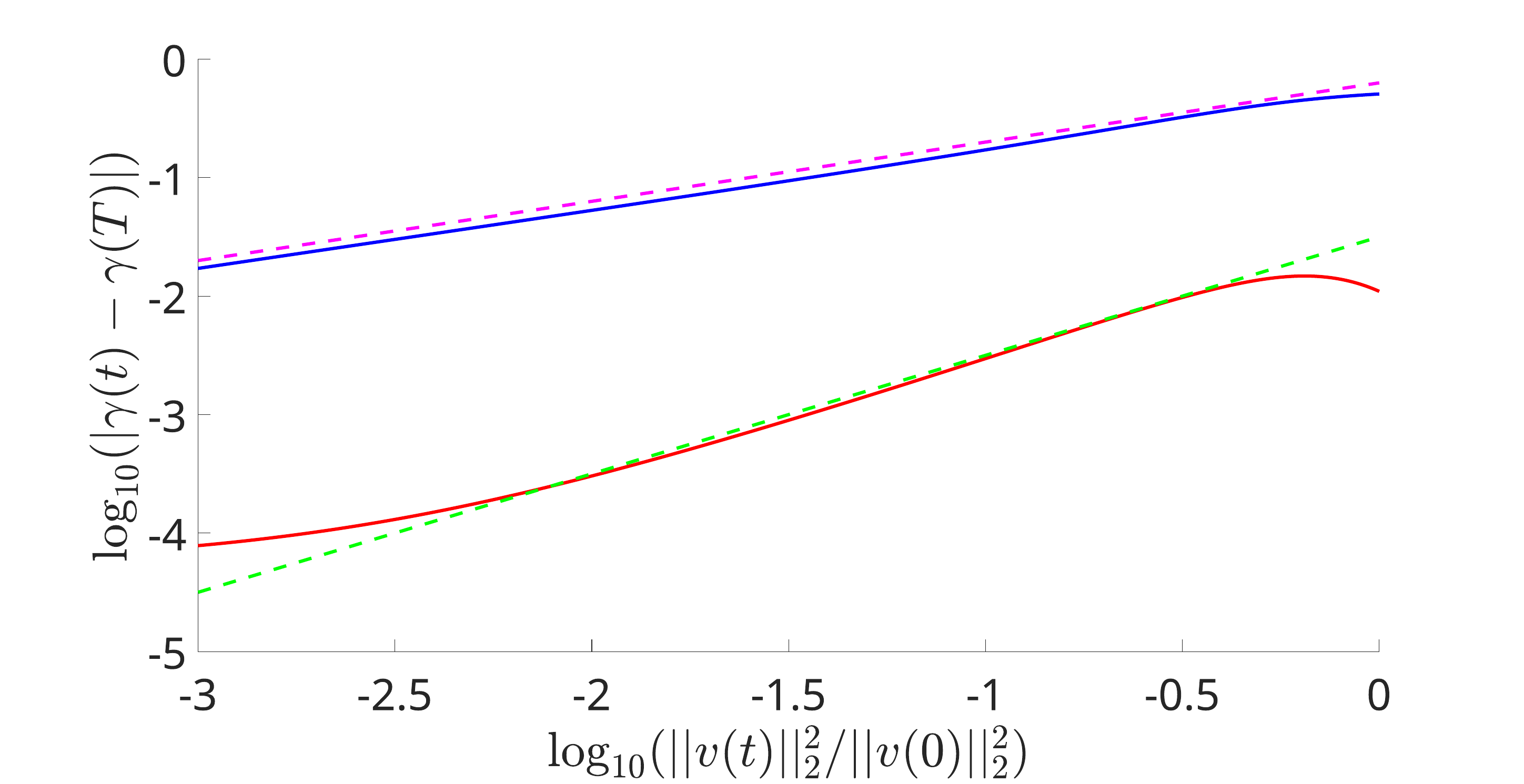}
    \caption{}
        \label{fig:DecayB}
\end{subfigure}
\caption{In both figures, the red line indicates the phase found by projecting on the adjoint eigenfunction $\psi$, and the blue line the phase found by projecting on the eigenfunction $\Phi'_0$. In Figure (a), the $L^2$-norm of the deviation $v(t)$ from the wave is is shown for both phases, showing a clear exponential decay. I normalized $\nrm{v(t)}_2^2$ by $\nrm{v(0)}_2^2$ in order to make both curves start at $\log_{10}(1)=0$. In Figure(b), the decay is not shown as a function of time, but as a function of the size of the phase $\gamma(t)-\gamma(50)$, where 50 is the endpoint of the simulation. The magenta dashed line has slope 1/2, while the green dashed line has slope 1.}
\label{fig:Decay}
\end{figure}

To drive this point home, let us conduct a numerical experiment. For both the projection on the eigenfunction and the adjoint eigenfunction, we compute $v(t)=u(\cdot-\gamma(t),t)-\Phi_0$. As initial condition, we take a small perturbation of $\Phi_0$. In Fig.~\ref{fig:Decay}, we show the norm of $v(t)$ for both phases, and indeed, the exponential decay is very similar. However, when we plot the decay of $\gamma(t)$ to $\gamma(\infty)$ as a function of $\nrm{v(t)}$, we see that slopes are very different for the two interpretations. Also note that $|\gamma(\infty)-\gamma(0)|$ is almost two orders of magnitude larger for the varational phase than for the H\&H phase.  

\subsubsection{Variational approach in weighted spaces}
We ended the previous section with the conclusion that the semigroup approach does not work, unless the linear operator is self-adjoint, which most examples are not. However, the linear operator can be forced to become self-adjoint by formulating the problem in a weighted space. In~\cite{stannatkruger}, which we will discuss in Section \ref{sec:phaselag}, the weight $\rho$ that makes the linear operator self-adjoint is found to be $\rho=\Phi_0'\psi^{-1}$. Hence, we have the identity
\begin{align*}
    \ip{v,\Phi'_0}_{L^2_\rho}=\ip{v,\psi}_{L^2},
\end{align*}
for any $v\in L^2_\rho(\R)$. Therefore, restricted to the smaller space $L^2_\rho(\R)$, the variational approach is identical to projecting directly on the adjoint eigenfunction. In many ways, this approach is similar to the approach by MacLaurin and Bressloff in~\cite{maclaurin2020}.

\subsubsection{Updated Inglis \& MacLaurin by MacLaurin.}
In~\cite{maclaurin2023phase}, MacLaurin solves the issues with the original paper by also projecting onto the adjoint eigenspace and explicitly including the constant $M\geq 1$ into the proofs. However, MacLaurin still refers to the phase as the variational phase, even though it does not solve the original variational problem (but then again, it would solve the variational problem in the smaller space $L^2_\rho$). Furthermore, MacLaurin states that the old and new variational phases are equivalent to leading order in $\nrm{u(t)-\Phi(\cdot-\gamma(t))}_2^2$. Given our analysis above and the results in Fig.~\ref{fig:Decay}, we now know this statement is false.\footnote{Upon request, MacLaurin specified that this statement just means ``the results are not affected by the choice of phase in any meaningful way," which I agree with.} 

The new phase is formulated in a more general setting than the old one. Specifically, the dimension of the symmetry group of the system is supposed to be $m$ instead of just 1, following the exposition in \cite{kapitula}. Hence, this allows for equations with more than 1 symmetry, such as spiral waves in 2d reaction-diffusion equations. In this new setting, and under a list of assumptions, MacLaurin shows stability of the waves on exponentially long time scales. However, it should be noted that MacLaurin does not explicitly show for all his examples that his assumptions (2.1)-(2.10) hold. 

All in all, in this new approach, the number of differences between MacLaurin and me has reduced to two: (i) I shift the solution, MacLaurin shifts the wave and (ii) I study deviations from $(\Phi_\s,c_\s)$, MacLaurin from $(\Phi_0,c_0)$. However, as noted in Remark 1, the second difference is not essential. 

\subsection{Physics literature}
The physics literature on itself is large enough to deserve its own summary; luckily, there is the invaluable book ``Noise in spatially extended systems" from '99, where an overview for reaction-diffusion equations is given. In~\cite[Sec. 6]{garciaspatiallyextended}, the following bistable reaction-diffusion equation is studied:
\begin{align}
    dU=[U_{xx}+f(U,a)]dt+\sigma g(U)\circ dW^Q_t.
\end{align}
Typically, we should think of an equation in the form of
\begin{align}
    dU=[U_{xx}+U(1-U)(U-a)]dt+\sigma U(1-U)\circ dW^Q_t.
\end{align}
Written in It\^o interpretation, this becomes\footnote{Note that the It\^o-Stratonovich correction term in the physics literature often lacks the $\frac{1}{2}$, because the spatial correlation is defined as $2q$ instead of $q$.}
\begin{align}
    dU=[U_{xx}+U(1-U)(U-a)+\frac{\s^2}{2}q(0)g'(u)g(u)]dt+\sigma u(1-u) dW^Q_t.
\end{align}
Here, I have to admit that I do not fully understand the computations that follow in~\cite{garciaspatiallyextended}, which is a nice way of saying that I believe the computations, as written down, are incorrect, e.g. Eq.~(6.18) makes no sense. However, I can follow the explanations around the equations, so what follows here is my interpretation of the calculations. First, we introduce a phase $\Gamma(t)=\Gamma_0+c_\mathrm{sne} t +\Delta(t)$, where $c_\mathrm{sne}$ is a yet undetermined speed and $\Delta$ is a stochastic process, and we shift to $\xi=x-\Gamma(t)$. I use \textit{sne} as an abbreviation for ``small noise expansion". In this frame, we assume that $U(\xi,t)=\Phi_\mathrm{sne}+\sigma U_1(t)+...$ and assume that $\Delta(t)=\sigma\Delta_1(t)+...$. When we collect terms of $\mathcal{O}(1)$ and the It\^o-Stratonovich correction, we find
\begin{align}    
\label{eq:TWEsne}
\Phi''_\mathrm{sne}+c_\mathrm{sne}\Phi'_\mathrm{sne}+f(\Phi_\mathrm{sne})+\frac{\s^2}{2}q(0)g'(\Phi_\mathrm{sne})g(\Phi_\mathrm{sne})=0.
\end{align}
Note that this is not the true lowest order. However, if you know the solution $(\Phi_0,c_0)$, it can be relatively straightforward to compute $(\Phi_\mathrm{sne},c_\mathrm{sne})$. 
Next, we find at $\mathcal{O}(\sigma)$
\begin{align}
    \mathcal{L}_\mathrm{sne} U_1=-\dot\Delta_1(t)\Phi'_\mathrm{sne}+g(\Phi_\mathrm{sne})dW^Q_t,
\end{align}
Here, $\mathcal{L}_\mathrm{sne}$ is the linearization of Eq.~\sref{eq:TWEsne} around $\Phi_\mathrm{sne}$ and it has a 0-eigenvalue with $\Phi'_\mathrm{sne}$ as eigenvector. Hence, the first-order equation for $U_1$ can only be solved when the rhs of the equation is orthogonal to the adjoint eigenfunction $\psi_\mathrm{sne}$, which still is given by $e^{c_\mathrm{sne}\xi}\Phi'_\mathrm{sne}$, up to normalization. Hence, we have to choose
\begin{align} 
\dot\Delta_1=\frac{\ip{\psi_\mathrm{sne},g(\Phi_\mathrm{sne})dW^Q_t}}{\ip{\psi_\mathrm{sne},\Phi'_\mathrm{sne}}}.
\end{align}
This expression can be used to compute the variance of the lowest-order perturbation of the phase. Upon replacing $\Phi_\mathrm{sne}$ with $\Phi_0$,  this indeed results in Eq.~\sref{eq:varphase}.

\subsubsection{Pros and cons}
The pros of this method are clear. You do not need any formal set-up, just some perturbation theory, and we have an approximation $(\Phi_\mathrm{sne},c_\mathrm{sne})$ for the average wave and an approximation for the variance in the movement of the front. The interesting part here is that the approximation $(\Phi_\mathrm{sne},c_\mathrm{sne})$ remains valid for all regimes, i.e. for $a\in [-1,1]$. The cons are, of course, that we do not know to what extent this approximation is valid. Or in fact, we do know, they are not valid at $\mathcal{O}(\sigma^2)$ (but they also never claim to be). The approximation $(\Phi_\mathrm{sne},c_\mathrm{sne})$ does match up nicely with the numerics, but this is mainly because it captures the drift induced by the Stratonovich interpretation. Any $\mathcal{O}(\sigma^2)$ interactions between the noise and the phase are not accounted for. This is a larger `problem' in the physics literature. Bifurcation diagrams such as Fig.~3 in~\cite{garcia2001} are drawn solely based on the effect of the It\^o-Stratonovich correction term. Why does this approach work? Well, the correction term $\frac{\s^2}{2}q(0)g'(\Phi_0)g(\Phi_0)$ is generally large, especially when we consider that the noise is often assumed to be white in space, i.e. $q(0)=\infty$. Of course, this does not make sense, and practically, the correlation length is chosen at the size of the spatial discretisation, e.g. $2q(0)=\frac{\s^2}{dx}$, for some small $\sigma$. This already hints at the fact that the Stratonovich SPDE cannot be understood in the space-time white noise limit without renormalization~\cite{Hairer2015}. 

\subsection{Cartwright and Gottwald: collective coordinate approach}
Cartwright and Gottwald developed a method to study stochastic bistable equations~\cite{cartwright2019} and stochastic KdV equations~\cite {cartwright2021collective}. In this so-called collective coordinate approach, the phase is not a priori set. 
A reference function $\hat u$ is introduced together with a set of coordinates $\mathbf{c}(t)$, and the Ansatz is made that $\hat u(x,t;\mathbf{c}(t))$ is a `good' approximation of the full solution. Good here means that the error $\mathcal{E}$ between the true solution and the Ansatz is small. For the SPDE
\begin{align}
    du=[u_{xx}+u(1-u)(u-a)]dt+u(1-u)d\beta_t,
\end{align}
where $\beta_t$ is a 1D Brownian motion, the Ansatz could be of the shape
\begin{align}
    \hat u(x,t;w(t),\phi(t))=\frac{1}{2}\left(1-\tanh(w(t)(x-\phi(t)))\right).
\end{align}
The function $w(t)$ represents the width of the profile and $\phi(t)$ the position, both solutions of an SDE with unknown coefficients. Now, the error is determined by choosing the coefficients of the SDEs for $w$ and $\phi$. The function $\hat u(x,t;\mathbf{c}(t))$ generates a 2D manifold indexed by the values of $\mathbf{c}(t)$ and we now assume that the dynamics of $\mathcal{E}$ is orthogonal to this manifold. In other words, $\mathcal{E}$ must be orthogonal to $\frac{\partial \hat u}{\p\phi}$ and $\frac{\partial \hat u}{\partial w}$. 

It is important to note here that this is not the same as the variational phase from~\cite{Inglis}, because we do not project onto $\Phi'_0(\cdot -\phi(t))$, but on $\frac{\partial \hat u}{\p\phi}$, which for the choice of $\hat u$ above reduces to
\begin{align}
    \frac{\partial \hat u(\cdot,t;w(t),\phi(t))}{\p\phi(t)}=\frac{1}{2}\text{sech}^2(w(t)(\cdot-\phi(t))).
\end{align}
When we demand that the deterministic and the stochastic part of $\mathcal{E}$ are orthogonal to both normal vectors, we find four equations for the four unknown coefficients of $\mathbf{c}(t)$. 
A stationary solution can be found of these equations that results in the same exact solution as found in~\cite{hamster2017}. 

Remark: The results of the collective coordinate approach are compared with numerical techniques, where the phase is defined as the point where the solutions attain the value $1/2$. 

\subsubsection{Pros and cons}
The main advantage of this approach is that you have a lot of flexibility. Note that $\hat u$ does not have to be the deterministic travelling wave. In principle, any function could work, but good results are only expected when $\hat u$ has the same shape as $\Phi_0$.  

A downside is the computational complexity of this approach. Where Hamster and Hupkes use the one symmetry of the Nagumo equation (translation), Cartwright and Gottwald introduce an extra equation for the width. For the KdV equation, Westdorp and Hupkes need 2 equations for the two symmetries of the system~\cite{westdorp2024long}, but Cartwright and Gottwald need 4. Furthermore, at this point, it is not clear how this approach could lead to any proofs, considering the validity of the approximations. 

Considering the noise, the equation is studied both with 1D Brownian motion and with a cylindrical process $W^Q_t$, where $q(x,y)\sim \exp(-|x-y|)$. However, this term is then multiplied by a function to localise the noise, i.e. to make the process trace-class, but this breaks the translation invariance. As far as I understand, this is not strictly necessary for the existence results on $\R$ and on the finite domain used for the numerical computations, translation-invariant kernels are of trace class anyway. 

\subsection{Kr\"uger, Stannat and Lang: Phase-lag method}
\label{sec:phaselag}
Stannat, Lang and Kr\"uger developed a technique that approximates a variational phase. 
Starting again with a deterministic PDE of the form
\begin{align}
    \frac{\p u}{\p t}=\partial_{xx}u+f(u),
\end{align}
in Inglis and Maclaurin~\cite{Inglis}, a phase $\gamma$ is introduced such that
\begin{align}
    \frac{d}{d\gamma(t)}\nrm{u(t)-\Phi_0(\cdot -\gamma(t))}^2=-2\ip{u(t)-\Phi_0(\cdot -\gamma(t)),\Phi'_0(\cdot -\gamma(t))}=0.
\end{align}
However, in a series of works by Stannat, Lang and Kr\"uger \cite{stannatnag,stannatbistable,stannatkruger,langstannat2016l2}, the following weaker version is introduced:\footnote{The order of explanation here is for educational purposes. In fact, this phase is older than the variational phase.}
\begin{align}
    \frac{d}{dt}\nrm{u(t)-\Phi_0(\cdot -\gamma(t))}=2\dot\gamma(t)\ip{u(t)-\Phi_0(\cdot -\gamma(t)),\Phi'_0(\cdot -\gamma(t))}\leq 0.
\end{align}
The inequality can be achieved by choosing 
\begin{align}
\label{eq:defphaselag}
    \dot\gamma_m(t)=-m\ip{u(t)-\Phi_0(\cdot -\gamma_m(t)),\Phi'_0(\cdot -\gamma_m(t))}
\end{align}
for any positive $m$. With this specific choice, the difference $\nrm{u(t)-\Phi_0(\cdot-\gamma_m(t))}$ must decay, given existence and regularity of the solution $(u(t),\gamma_m(t))$. We did not specify a specific Hilbert space yet, and typically $L^2$ or a weighted version $L^2_\rho$ is used. 

An advantage of this approach is that there is no $\ip{\Phi_0'+v',\psi}^{-1}$-term as there is Eq.~\sref{eq:defav}. This implies that $\gamma_m(t)$ circumvents the possibilities of singularities in $\dot\gamma$, which corresponds to discontinuities in the phase description. In the case of the Nagumo equation, this property can be really useful as the numerics indicates that stochastic travelling waves persist after a discontinuity in the phase. However, in the case of the FitzHugh-Nagumo equation, stochastic perturbations might destroy the wave back to the background state, and no phase description will recover from this. 

In~\cite{stannatnag,stannatbistable}, Stannat proves that this approach works for deterministic bistable reaction-diffusion equations. The proof depends crucially on the following estimate, referred to as a `local dissipativity estimate':
\begin{align}
    \ip{u_{xx}+f'(\Phi)u,u}_{L^2(\R)}\leq k_*\nrm{u}_{H^1}+C_*\ip{u,\Phi'}_{L^2(\R)}^2,
\end{align}
for some specifically chosen positive $k_*$ and $C_*$. 
In~\cite{stannatbistable}, this estimate is proved using extensive computations. Currently, I am not aware of similar results in a more abstract setting. The power of this inequality is that, once it has been shown, the proof of the orbital stability of the travelling wave is straightforward.

In \cite{stannatkruger}, for the SPDE
\begin{align}
    dU=[U_{xx}+f(U)]dt+\varepsilon dW^Q_t,
\end{align}
the following decomposition\footnote{Note that in most works by Stannat et al., the role of $v$ and $u$ is switched.} is made:
\begin{align}
    V(t)=U(t)-\Phi_0(\cdot-c_0t-C_m(t)),
\end{align}
where $C_m(t)$ solves an equation like \sref{eq:defphaselag}. The noise is additive and $Q$ is of trace class. It is important to stress that phase $C_m(t)$ solves an RODE, and hence the speed $c_m(t):=\dot C_m(t)$ is a well-defined stochastic process, which solves an SDE. When we apply the chain rule to the equation for $V(t)$ above, we only take derivatives of $C_m(t)$ and therefore no It\^o calculus is needed.
Next, equations for $V(t)$ and $C_m(t)$ are computed and it is shown that these equations can be meaningfully expanded in $\e$ up to certain timescales. All analysis is done in a weighted $L^2$ space with weight $\rho(x)=e^{c_0x}$. Hence, everywhere where the authors write $\ip{v,\Phi'_0}_\rho$, we can read $\ip{v,\psi}$. Next, the authors show that, for $m\to\infty$, $C_m(t)$ converges to a projection onto the adjoint eigenspace, and in turn show that this projection onto the eigenspace is a `good' approximation of the phase on short timescales. The authors also formally recover Eq.~\sref{eq:varphase} provided that $Q$ is approximately translation invariant. The translation-invariant case itself is not covered, since they assume that $Q$ is of trace class.

\subsubsection{A comment on Stannat's original papers}
In~\cite{stannatnag,stannatbistable}, the deterministic proof is extended to the stochastic case. However, the result claims that deviations from the stochastic travelling wave remain small infinitely long with nonzero probability. This result cannot be true, as the same result also implies that a solution can escape on the finite time interval $[0,T]$ with nonzero probability. The problem stems from the fact that (in the notation from~\cite{stannatbistable}) for the martingale
\begin{align}
    M_t=\int_0^t\ip{v(t),g(v(t)+\Phi(\cdot +ct+C(t))dW^Q_s},
\end{align}
and stopping time 
\begin{align}
    \tau:=\inf\{t\geq 0 | \nrm{v(t)}\geq c_*\},
\end{align}
we cannot apply Doob's optimal stopping theorem to the integral
\begin{align*}
    I(t\wedge\tau)=\int_0^{t\wedge\tau}e^{-k_*(t\wedge\tau-s)}dM_s,
\end{align*}
as it is not a local martingale. Hence, $E[I(t\wedge\tau)]$ can be nonzero, while it is assumed to be zero in~\cite{stannatnag,stannatbistable}. 
In later works, such as~\cite{stannatkruger,lang,langstannat2016l2}, but also~\cite{eichinger2022multiscale,gnann2025solitary}, martingales such as $I$ above are bounded by bounding $E[\sup_{0\leq t\leq T}I(t)]$ using e.g. a Burkholder-Davis-Gundy inequality. This circumvents the problems with Doob's optimal stopping theorem. 
    
\subsubsection{Pros and Cons}
For now, most computations are relatively straightforward, as the noise is additive in \cite{stannatkruger,eichinger2022multiscale}. The phase is described by a random ODE instead of an SDE, which also reduces the computational complexity. Furthermore, in~\cite{eichinger2022multiscale}, the authors do not need an analytic semigroup, a $C_0$-semigroup suffices, which implies that they can treat an FHN equation without diffusion in the second component. 

A disadvantage is that there are now two perturbation parameters, $\sigma$ and the relaxation rate $m$. It is not a priori clear how we should approach a systematic expansion around $\sigma$ and $m^{-1}$. Furthermore, all bounds at this point are polynomial in time instead of exponential. It is not clear at this point to what extent an approximate phase like $C_m$ can lead to bounds at exponentially long timescales. All the proofs at this point explicitly depend on the fact that $Q$ is trace-class. Again, it is not immediately clear if the proofs can be extended to translation-invariant $Q$. In the case of multiplicative noise, this problem might be solved more easily. Note that the technique has recently been extended to infinite-dimensional noise with H\"older-continuous paths by Saef and Stannat~\cite{saef2025stability}, where translation-invariance does not seem to be an issue.

\section{Miscellaneous}

\subsection{K\"uhn: Integral}
In~\cite{kuehnreview}, K\"uhn introduces an integral to study the location of a travelling front on a bounded domain. The idea is simple: when the front moves, the area under the curve increases or decreases with the movement. Hence, a change in the area under the curve indicates front movement. Numerically, this is very straightforward to implement, once you have the SPDE simulation. However, this approach only works on a finite domain for finite timescales, i.e. until you hit the boundary. Note that on a finite domain $D$, the integral can be understood as a functional on $L^2(D)$ and hence is, at least deterministically, equivalent to all other phases described by an $L^2(D)-$functional. 

\subsection{Mueller et al: Level sets}
Muellet et al.~\cite{mueller1995, mueller2011} developed in a series of papers a technique to prove the Brunot-Derrida conjecture. Therefore, let us write down a F-KPP equation with square root noise:
\begin{align}
    du=[u_{xx}+u(1-u)]dt+\sigma\sqrt{u(1-u)}dW^Q_t.
\end{align}
For this noise term, the solutions have the property that they are `compact' in the sense that they converge to 0 and 1 in finite space~\cite{shiga1994two}. Hence, we can write down an interval $[a(t),b(t)]$ where the travelling wave decays from 1 to 0, and these values are finite. This choice allows the authors to define a stochastic speed
\begin{align}
    c_\sigma=\lim_{t\to\infty}t^{-1}b(t).
\end{align}
Note that there is no expectation, i.e. the limit of the random process is not random. Furthermore, it is shown that the speed $c_\sigma$ can be expanded as
\begin{align}
    c_\sigma=c_0+\frac{\pi^2}{|\log(\sigma^2)|^2}+\mathcal{O}\left(\frac{\log(\log(\sigma))}{|\log(\sigma^2)|^3}\right).
\end{align}
Interestingly, none of the other approaches mentioned in this overview has been used to prove the validity of this expansion. Note that this does not mean that it is impossible, it just highlights that proving the validity of this expansion is a difficult problem.   

\subsection{Reverse phase tracking}
We can systematically look for exact solutions of SPDEs by following an approach that I will call `reverse phase tracking'. Reverse in the sense that we do not find an Ansatz for an SPDE, but rather start with an Ansatz and find an SPDE. 

\subsubsection{Advective 1D noise}
Suppose we have the following general PDE forced by a 1d Brownian motion in Stratonovich sense:
\begin{align}
    \label{eq:ExactSolGeneralS}
    dU=[F(U,U_x,U_{xx},...)]dt+\sigma U_x\circ d\beta_t,
\end{align}
and introduce the phase
\begin{align}
    d\Gamma(t)=cdt+\sigma d\beta_t.
\end{align}
Then, $v(t)=U(\cdot-\Gamma(t),t)$ solves the PDE
\begin{align}
    v_t=F(v,v_x,v_{xx},..)+cv_x.
\end{align}
This follows directly from the classic chain rule because the noise is interpreted in the Stratonovich sense. Therefore, if the deterministic version of Eq.~\sref{eq:ExactSolGeneralS} has a travelling wave, we can immediately solve the stochastic version. When we look at the It\^o interpretation, we introduce a phase 
\begin{align}
    d\Gamma(t)=c_\s dt+\sigma d\beta_t
\end{align}
and find that $v(t)=U(\cdot-\Gamma(t),t)$ solves the PDE
\begin{align}
    v_t=F(v,v_x,v_{xx},..)+c_\sigma v_x-\frac{\s^2}{2}v_{xx}.
\end{align}
Hence, we can exactly solve the It\^o version of the SPDE when we can solve the perturbed TWE equation above. For example, we can solve the following stochastic Nagumo equation:
\begin{align}
    dU=[dU_{xx}+U(1-U)(U-a)]dt+\sigma U_xd\beta_t,
\end{align}
because we now know that $\Phi_\s(\cdot-\Gamma(t),t)$ is a solution when
\begin{align}
    (d-\frac{1}{2}\sigma^2)\Phi''_{\s}+c_\sigma\Phi'_\s+\Phi_\s(1-\Phi_\s)(\Phi_\s-a)=0.
\end{align}
Solving this equation is straightforward when $d-\frac{1}{2}\s^2>0$. Note that the boundary conditions do not change, as $U_x$ vanishes at the deterministic endpoints 0 and 1. However, if we were to use this approach for the KdV equation, as was done in \cite{cartwright2021collective}, we would solely get a negative diffusive term. We cannot solve the resulting TWE exactly, but we still conclude that the soliton is unstable due to the negative diffusive term. 

In the case of the following stochastic Burgers-KdV equation
\begin{align}
    dU=[\delta U_{xxx}+\mu U_{xx}+\beta UU_x]dt+\sigma U_xd\beta_t,
\end{align}
we can solve this SPDE by solving the following perturbed TWE:
\begin{align}
    \delta \Phi'''_\s+(\mu-\frac{\s^2}{2}) \Phi''_\s+c_\s\Phi'_\s+\beta \Phi_\s\Phi_\s'=0.
\end{align}
This time, the equation is stable and solvable as long as $(\mu-\frac{\s^2}{2})>0$. 
Note that for this specific equation, the speed $c$ is uniquely determined by the endpoints of the wave $\Phi_\s$:
\begin{align}
    c_\s=-\frac{\beta}{2}(\Phi_\s(-\infty)+\Phi_\s(\infty)).
\end{align}
This can be checked by integrating over $\mathbb{R}$ and assuming that all derivatives vanish at $\pm\infty$. Hence, $c_\s=c_0$ when we keep the endpoints fixed, but the profile $\Phi_\s$ is changed. However, in \cite{adjibi2024exact}, this is not accounted for. There, the distinction between $\mu$ and $\mu_{\rm eff}=\mu-\s^2/2$ is not made, and the boundary conditions are not fixed. For the deterministic travelling wave presented in \cite{adjibi2024exact}, the boundary conditions depend on $\mu$. Hence, when we want to understand how the noise changes the solutions without changing the boundary conditions, I'm not sure we can find an exact solution.

\subsubsection{Additive 1D noise}
Let us study the following Burgers-KdV equation:
\begin{align}
    dU=[\delta U_{xxx}+\mu U_{xx}+\beta UU_x+\alpha U_x]dt+\sigma(t)dW_t.
\end{align}
Following the computations in~\cite{adjibi2024exact}, we choose 
\begin{align}
    U(t)=\Phi_\s\left(\cdot-c_\sigma t+\beta\int_0^t\int_0^s\sigma(s')dW_s'ds\right)+\int_0^t\sigma(s)dW_s.
\end{align}
This results in
\begin{align}
    \left[-c_\s+\beta\int_0^t\sigma(s)dW_s\right]\Phi'_\s=\delta\Phi'''_{\s}+\mu \Phi''_{\s}+\beta \left(\Phi_\s+\int_0^t\sigma(s)dW_s\right)\Phi'_\s,
\end{align}
which reduces to
\begin{align}
    -c_\s\Phi'_\s=\delta\Phi'''_{\s}+\mu\Phi''_{\s}+\beta\Phi_\s\Phi'_\s.
\end{align}
Hence, the TWE is independent of $\sigma(t)$, and the shifted deterministic wave $(\Phi_0,c_0)$ plus Brownian motion solves the SPDE. Note that this method works because of the structure of the PDE. When we would study a modified KdV-Burgers equation, e.g. by replacing $uu_x$ with $u^2u_x$, this Ansatz wouldn't work anymore, but you can add any type of nonlinearity depending on spatial derivatives.

\section{Isochronal Phase}
Here, we study the isochronal phase and some related results. Especially, in Section \ref{sec:Comparison}, we will study the relation between the results of Hamster and Hupkes, and the isochronal phase. Next, we will discuss Van Winden's predicted phase, which is related to but not identical to the isochronal phase. However, I think that the predicted phase will feel natural after the discussion in Section \ref{sec:Comparison}.

\subsection{Adams \& MacLaurin: Isochronal Phase}
Suppose we have a deterministic reaction-diffusion equation, where orbital stability of the travelling wave has been established. This means that for any function $u_0$ close enough to the manifold $\mathcal{M}$ of travelling waves, we know that there exists a number $\gamma_\infty(u_0)$ such that, in the travelling wave coordinate $\xi=x-c_0t$, a solution $u(t)$ of the equation with $u_0$ as initial condition converges to $\Phi(\cdot+\gamma_\infty)$~\cite{kapitula}. Hence, the map $\pi$ that maps $u_0$ to $\gamma_\infty(u_0)$ is well-defined for all $u_0$ in the basin of attraction of the manifold $\mathcal{M}$. This map $\pi$ is called the isochron map, and defines the isochronal phase. By definition, the isochronal phase does not describe where you are now, but where you are at $t\to\infty$. An important property of the Isochronal phase is that it is invariant under the flow of the PDE. When we define $\phi_t(u_0)$ as the flow of the PDE starting at $u_0$, then by definition of the isochronal phase, it is immediately clear that $\pi(\phi_t(u_0))=\pi(u_0)$ for all $t$. 

As a historical side note, the isochronal phase was originally designed to study oscillators, see~\cite{winfree1974patterns,guckenheimer1975isochrons}. Here, the isochronal phase allows you to make a very natural connection between a point moving exactly along the periodic orbit and a point near the periodic orbit (assuming the periodic orbit is stable). This approach is also naturally extended to the stochastic case, see the results by Adams~\cite{adams2023asymptotic}.   

Here, we focus on the application of the isochronal phase to travelling waves in stochastic PDEs. We can apply the isochronal phase directly to the stochastic version of the PDE, because the map $\pi(U_t)$ is well defined as long as the stochastic solution remains within the (deterministic) basin of attraction of the manifold $\mathcal{M}$. Note, however, that $\pi$ is not invariant under the stochastic flow. Even though $\pi$ is straightforward to define for the stochastic solutions, computing it still requires us to solve the deterministic PDE for each point in time. 

What can we now do with this phase? We start by stating the SPDE
\begin{align}
    dX_t=[AX_t+N(X_t)]dt+\sigma B(X_t)dW_t,
\end{align}
either on $\R$ or on a periodic domain, large enough for the existence of e.g. pulses of the FHN-equation to exist. Adams and MacLaurin~\cite{adams2024existence,adams2025isochronal} show that you can write down an It\^o equation for $\pi$:
\begin{align}
\begin{split}
    \pi(X_{t\wedge\tau})=\,&\pi(X_0)+\int_0^{t\wedge\tau}D\pi(X_s)[A X_s+N(X_s)]ds+\frac{\sigma^2}{2}\int_0^{t\wedge\tau}\sum_{k\in\mathbb{N}}D^2\pi(X_s)[B(X_s)e_k,B(X_s)e_k]ds\\
    &\quad+\sigma\int_0^{t\wedge\tau}D\pi(X_s)B(X_s)dW_s,
\end{split}
\end{align}
where $\tau$ is the stopping time of leaving the basin of attraction. This formula is exactly what you expect, but the authors show that the first and second Fr\'echet derivatives actually exist, and the stochastic integral is well-defined. Note that this equation holds in the strong sense, which is anything but trivial, especially as the SPDE itself does not have solutions in the strong sense.  

When we now focus on an equation of FHN-type on a periodic domain, we can understand the phase $\pi(X_{t\wedge\tau})$ as a stochastic process on the circle. Note that for the following results by Adams and MacLaurin, it is essential to work on a compact domain like the circle. The question we can ask now is whether or not we can describe the movement of the pulse at intermediate timescales. At short timescales, the position of the wave is very much determined by the (deterministic) initial condition. However, we cannot take the limit of $t\to\infty$ as solutions of the FHN equation will stop being soliton-like at some point with probability 1.  If we find a stationary distribution that accurately describes the distribution of $\pi_t$ at intermediate timescales and conditioned on not leaving the basin, we call it a quasi-stationary distribution (QSD) of $\pi_t$. See \cite{adams2025isochronal} for a rigorous definition.  

Directly computing such a QSD is hard, hence in \cite{adams2025isochronal}, see Theorem 4.2, the procedure is split into two. First, we write down an approximation $\bar\pi_t$ of $\pi_t$ for which it is easier to show that it has a stationary measure. Next, it is shown that $\pi_t$ indeed remains close to this stationary measure at intermediate timescales. 

To find an approximation $\bar\pi_t$ for $\pi_t$, we note that at this point the isochronal phase is dependent on the full solution $X_t$. Hence, the goal is to find an approximation of $\pi_t$ that is independent of $X_t$. This can be achieved by assuming that the noise is small, and therefore the solution $X_t$ is close to the manifold, which allows us to approximate $X_t\approx\Phi_0(\cdot-\pi(X_t))$. Hence, we can now define a process $\bar\pi_t$ that solves the SDE
\begin{align}
\begin{split}
    \bar\pi_t=\,&\bar\pi_0+\underbrace{\int_0^tD\bar\pi_s[A\Phi_0(\cdot-\bar\pi_s)+N(\Phi_0(\cdot-\bar\pi_s))]ds}_{=0}+\frac{\s^2}{2}\int_0^t\sum_{k\in\mathbb{N}}D^2\bar\pi_s[B(\Phi_0(\cdot-\bar\pi_s))e_k,B(\Phi_0(\cdot-\bar\pi_s))e_k]ds\\
    &\quad+\sigma\int_0^tD\bar\pi_sB(\Phi_0(\cdot-\bar\pi_s))dW_s.
\end{split}
\end{align}
Adams and MacLaurin do not name this SDE, but in \cite{kuehn2025synchronization}, this SDE is called the autonomous isochronal phase. Note that the first integral disappears because we work in the travelling wave coordinate, i.e. $A=\partial_{xx}+c_0\partial_x$. It can now be shown that this process has an invariant measure and that the full process $\pi(X_t)$ remains, after a short initial dynamics, close to this stationary measure on exponentially long timescales. 

\subsubsection{Heuristics for the QSD}
Without following the computations by Adams and MacLaurin, can we get an intuition of what the stationary distribution must be? We know that, at least on $\R$, any reasonable phase description is described by a scaled Brownian motion at lowest order, as discussed before. Hence, the phase has no stationary distribution on $\R$. However, projected onto a compact set like the circle, Brownian motion does have a stationary distribution, and it is the uniform distribution. Therefore, we expect the QSD to be uniform as well. 

To take an even more general picture: we know that the SPDE we are studying is translation invariant. Hence, we expect the QSD for the phase, if it exists, to be also translation invariant, which in turn implies that the QSD is uniform. 

\subsubsection{Pros and cons}
The power of the isochronal phase becomes clear in the following result, formulated in words: \textit{As long as the solution remains inside the basin of attraction, the dynamics can be well-approximated (after an initial transient) by an invariant measure}. Furthermore, this approach can be extended to Banach spaces, which opens up our toolbox to a whole range of new equations, especially as the phase description and the solution theory for the SPDE do not have to be formulated in the same space~\cite{kuehn2025synchronization}. 

The problem is that up to the writing of this document, there was no clear path to compute the QSD of the isochronal phase. First, we would need to find an equation that matches all the technical conditions in Adams. I assume the FHN equation satisfies these conditions, especially the doubly diffusive version, but this is not shown explicitly. In Section \ref{sec:Comparison} below, we show that the lowest-order expansions such as those obtained by Hamster \& Hupkes for the Nagumo equation, overlap with the autonomous isochronal phase when posed on the real line.  

Lastly, the results in \cite{adams2025isochronal} explicitly depend on the fact that we a priori know that the solutions remain in the basin of attraction for an exponentially long time. This is shown by applying techniques from \cite{maclaurin2023phase}. I assume such a result could also be shown directly using the isochronal phase, but this is not found in the literature. 

\subsection{Comparison of stochastic freezing and the autonomous isochronal phase}
\label{sec:Comparison}
For simplicity, we will study a one-component SPDE like the following Nagumo equation forced by an infinite-dimensional process:
\begin{align}
dU &= \big[ U_{xx} + f_{\mathrm{cub}}(U)]dt
+ \sigma g(U) d W^Q_t.
\end{align}
We assume $g(0)=g(1)=0$, as we wish to remain within the H\&H framework, but for the formal computations here, this will not really matter. Suppose we have a state $\bar U$ of the SPDE, what is the isochronal phase $\pi(\bar U)$? That is the value $\gamma_\infty$ where $\Phi_0(\cdot-c_0t-\gamma_\infty)$ is the limit of $t\to\infty$ for the deterministic flow of the FHN equation with $\bar U$ as initial condition.\footnote{That is, we use $\Phi_0$ as reference function for the phase 0 and assume that $\pi(\Phi_0)=0$. If you choose a reference function with $\pi(\Phi_{\rm ref})\neq0$, you'll have to think in terms of `changes in phase'.} We solve the PDE with $\bar U$ as initial condition using the approach from Section \ref{sec:HH}, i.e. we solve
\begin{align}
    \begin{split}
        v(t)=&S(t)v_0+\int_0^tS(t-s)[N(v(s))+a(v(s))\p_x(\Phi_0+v(s))]ds\\
        \gamma(t)=&\gamma_0+\int_0^ta(v(s))ds
    \end{split}
\end{align}
However, we have to be precise here. Let us therefore explicitly write down the relation between $\bar U$ and the pair $(v_0,\gamma_0)$: We define $\gamma_0$ such that 
\begin{align}
\label{eq:comp:gamma0}
    \ip{\bar U(\cdot+\gamma_0)-\Phi_0,\psi}=0,
\end{align}
and once we have solved $\gamma_0$, we set 
\begin{align}
    v_0=\bar U(\cdot +\gamma_0)-\Phi_0.
\end{align}
Note that we showed in \cite{hamster2017} that $\gamma_0$ exists and is unique when $\bar U$ is close to $\Phi_0$. However, it should be straightforward to show that $\gamma_0$ exists when $\bar U$ is close to any translate of $\Phi_0$.  
Now, given the fact that the limit of the deterministic flow is unique, we must have that
\begin{align}
\label{eq:DefIphaseHH}
    \pi(\bar U)=\gamma_\infty=\gamma_0+\int_0^\infty a(v(t))dt.
\end{align}
Following MacLaurin and Adams, we know that $\pi$ is a functional with smooth Fr\'echet derivatives, and we have the following It\^o formula:
\begin{align}
\begin{split}
    \pi(U(t))=&\pi(U_0)+\underbrace{\int_0^t\pi'(U(s))[\p_{\xi\xi}U(s)+c_0\p_{\xi}U(s)+f(U(s))]ds}_{=0?}+\int_0^t\pi'(U(s))[g(U(s))dW^Q_s]\\
    &+\frac{\sigma^2}{2}\sum_{i=0}^\infty\int_0^t\pi''(U(s))[g(U(s))\sqrt{Q}e_i,g(U(s))\sqrt{Q}e_i]ds.
    \end{split}
\end{align}
This version comes from \cite{dapratomild}, but note that the term with the underbrace is zero in the version of \cite{adams2025isochronal}, which in turn is used by \cite{kuehn2025synchronization}. I don't know why this is. 

At this point, we haven't achieved much, as knowing $\pi(U(t))$ is as complicated as knowing $U(t)$. Therefore, the approach taken in \cite{adams2024existence,kuehn2025synchronization} is to define an approximation $\pi_*(t)$ of $\pi(t)$ where we replace $U(t)$ with the approximation $\Phi_0(\cdot-\pi_*(t)):=\Phi_0(t)$:
\begin{align}
\begin{split}
    \pi_*(t)=&\pi_0+\underbrace{\int_0^t\pi'(\Phi_0(s))[\p_{\xi\xi}\Phi_0(s)+c_0\p_{\xi}\Phi_0(s)+f(\Phi_0(s))]ds}_{=0}+\s\int_0^t\pi'(\Phi_0(s))[g(\Phi_0(s))dW^Q_s]\\
    &+\frac{\s^2}{2}\sum_{i=0}^\infty\int_0^t\pi''(\Phi_0(s))[g(\Phi_0(s))\sqrt{Q}e_i,g(\Phi_0(s))\sqrt{Q}e_i]ds.
    \end{split}
\end{align}
This time, we do see that the first integral disappears, as $\p_{\xi\xi}\Phi_0+c_0\p_{\xi}\Phi_0+f(\Phi_0)=0$ for all translates of $\Phi_0$.

Our first task is to compute the Fr\'echet derivative $\pi'(\Phi_0(s))[w]$. This has been done in \cite[Lem. 3.5]{kuehn2025synchronization} and it turns out that
\begin{align}
\label{eq:comp:FirstFder}
    \pi'(\Phi_0)[w]=-\ip{w,\psi}.
\end{align} When we take a step back and rephrase the question ``what is the first Fr\'echet derivative of the isochronal phase" as ``suppose I take a step $w$ away from $\Phi_0$, what is the limit of the linearised dynamics", we can find the answer already in \cite{Sattinger} from '76.  
\begin{remark}
    NB: The first Fr\'echet derivative in \cite{kuehn2025synchronization} does not have the minus sign. Is this a computation error or a matter of definition? However likely a computation error might be, I believe the sign difference comes from a sign choice. In \cite{kuehn2025synchronization}, the projection onto the zero eigenspace is defined as $Pf=-\ip{f,\psi}\Phi'_0$. Hence, their adjoint is normalized such that $\ip{\psi,\Phi'_0}=-1$ instead of $\ip{\psi,\Phi_0'}=1.$
\end{remark}

Now, at order $\sigma$, we find
\begin{align}
\begin{split}
    \pi_*(t)=-\sigma\int_0^t\ip{g(\Phi_0(s))dW^Q_s,\psi(\cdot-\pi_*(s))}&=-\sigma\int_0^t\ip{g(\Phi_0)dW^Q_s(\cdot+\pi_*(s)),\psi}\\
    &\approx-\sigma\int_0^t\ip{g(\Phi_0)dW^Q_s,\psi},
\end{split}
\end{align}
where the approximation holds on timescales where the realization $W_t^Q$ and its shifted counterpart $W^Q_t(\cdot+\pi_*(t))$ remain close. However, note that the statistical properties of this integral are independent of the specific shift.

As I will show in Appendix~\ref{app:Iso}, the same expression for the first Fr\'echet derivative can be found by computing the Fr\'echet derivative of Eq.~\sref{eq:DefIphaseHH}. K\"uhn and Van Winden do not report the second Fr\'echet derivative. Again in Appendix~\ref{app:Iso}, I explicitly compute the second Fr\'echet derivative. It is important to note two different contributions, coming from $\gamma_0$ and from $\int_0^\infty a(v(t))dt$. First, we find that
\begin{align}
    D^2\gamma(\Phi_0(s))[w,w]=-\ip{\Phi_0'',\psi}\ip{w,\psi(s)}^2-2\ip{w,\psi'(s)}\ip{w,\psi(s)}.
\end{align}
Plugging this back into the It\^o formula results in
\begin{align*}
    \frac12\sum_{i=0}^\infty D^2\gamma(\Phi_0(s))&[g(\Phi(s))\sqrt{Q}e_i,g(\Phi(s))\sqrt{Q}e_i]\\&=-\frac12\sum_{i=0}^\infty\ip{\Phi_0'',\psi}\ip{g(\Phi(s))\sqrt{Q}e_i,\psi(s)}^2-\ip{g(\Phi(s))\sqrt{Q}e_i,\psi'(s)}\ip{g(\Phi(s))\sqrt{Q}e_i,\psi(s)}\\
    &=-\frac12\ip{\Phi_0'',\psi}\ip{Qg(\Phi_0)\psi,g(\Phi_0)\psi}-\ip{Qg(\Phi_0)\psi',g(\Phi_0)\psi}.
\end{align*}
The details of the second step are found in the Appendix. The result is identical to the term $c_{0;2}$ from Hamster and Hupkes, i.e. the second-order correction to their instantaneous wave speed, see Eq.~(3.14) in \cite{Hamster2020}.  

Next, we study the integral $\int_0^\infty a(v(t))dt:=J(\bar U)$, Note that $a$ is nonlinear and hence $DJ(\Phi_0(s))[w]=0$. The second order derivative is found to be
\begin{align}
    D^2J(\Phi_0(s))[w,w]=\int_0^\infty\ip{f''(\Phi_0)\left(S(s')[w(\cdot+\pi_*(s))-\Phi_0'\ip{w,\psi(\cdot-\pi_*(s))}]\right)^2,\psi}ds'.
\end{align}
Upon plugging in $w=g(\Phi_0(s))\sqrt{Q}e_i$ and taking the infinite sum, we get from the It\^o formula
\begin{align}
    \frac12\sum_{i=0}^\infty\int_0^\infty\ip{f''(\Phi_0)\left(S(s)[g(\Phi_0)\sqrt{Q}e_i-\Phi_0'\ip{g(\Phi_0)\sqrt{Q}e_i,\psi}]\right)^2,\psi}ds.
\end{align}
This is exactly the second-order correction to the orbital drift from Hamster and Hupkes. Hence, we conclude that, on average, the autonomous isochronal phase and the phase-tracking approach overlap at second order.

\subsection{Van Winden: Predicted phase}
In \cite{vanWinden2024noncommutative}, the predicted phase is formulated in a quite abstract setting of Banach spaces and Lie groups. Before we go there, let's formulate the concept for people who feel uneasy when geometry has left the room. From the discussion in the previous section on the first Fr\'echet derivative of the isochronal phase, we can conclude the following result for the deterministic case. Suppose we have as initial condition $u_0:=\Phi_0+v_0$ for $v_0$ small. Then, after some time $t$, the solution $u(t)$ will be close to $\Phi_0(\cdot-c_0t-\ip{v_0,\psi})$. Van Winden now defines the predicted phase $\gamma(t)$ as
\begin{align}
    \gamma(t)=c_0t+\ip{v_0,\psi}.
\end{align}
Note that this phase is not exact, but the power lies in the fact that it can be computed explicitly (in fact, we did, it is the equation above). Furthermore, we know that on a sufficiently large interval $[0,T]$ the difference between $u(T)$ and $\Phi_0(\cdot-\gamma(T))$ is reduced from $\mathcal{O}(\nrm{v_0})$ to $\mathcal{O}(\nrm{v_0}^2)$. This means that the error is reduced again on $[T,2T]$, $[2T,3T]$, etc. Hence, the stability proof can be closed by a nonlinear iteration argument. 

Now, for the Nagumo equation, orbital stability is of course not new. However, in \cite{vanWinden2024noncommutative}, the proof is formulated in a much more general setting. By formulating the problem in Banach spaces and treating the symmetries of the system (which inherently lead to orbital stability instead of `normal' stability) in terms of (non-commutative) Lie groups, a large class of deterministic equations can be studied at once. Furthermore, note that the semigroup generated by the linear operator is assumed to be $C_0$ instead of analytic, allowing for the study of single-diffusive FHN-equations. 

In the stochastic case, the solutions are studied in 2-smooth Banach spaces.\footnote{No, I do not claim to understand what this means. According to \cite{vanWinden2024noncommutative}, $L^p$-spaces have this property for $p\in[2,\infty)$.} This is a little more restrictive than the deterministic case, but still far more general than all versions discussed here that rely on Hilbert spaces. How is the predict phase used in the stochastic case? Suppose we want to know something about the solution $U(t)$ at $t=(n+1)T$. Then, using the predicted phase, we can estimate properties of $U((n+1)T)$ using information on $U(nT)$. In turn, we can estimate $U(nT)$ by understanding $U((n-1)T)$ etc. This leads to bounds similar to the ones in \cite{hamster2020exp}. However, it seems that the proof is significantly shorter. \\

\noindent\textbf{Competing Interests}\\
There are no competing interests. \\

\noindent\textbf{Acknowledgements}\\
The work of CH is funded by the Dutch Institute for Emergent Phenomena (DIEP) at the University of Amsterdam via the program Foundations and Applications of Emergence (FAEME). CH would like to thank Zach Adams for the input on this document.

\bibliographystyle{klunumHJ}
\bibliography{ref} 

\appendix
\section{Computation of the first and second Fr\'echet derivatives. }
\label{app:Iso}
In the main text, we showed the Fr\'echet derivatives of Eq.~\sref{eq:DefIphaseHH}. Note that we will differentiate with respect to the deviation from the translate $\Phi_0(s)$. We now write $\bar U=\Phi_0(s)+w$, i.e. we assume that $\bar U$ is close to the translate $\Phi_0(\cdot-\pi_*(s))$ as $w$ is small. Note however, that this does not mean that $\gamma_0(\bar U)$ is small, it is close to $\pi_*(s)$. Specifically: $\gamma_0(\Phi_0(s))=\pi_*(s)$.

\subsection{Derivatives of $\gamma_0$}
As a first step, we shift $\bar U$ with $\gamma_0(\bar U)$ to disconnect them: 
\begin{align*}
        \ip{\bar U(\cdot+\gamma_0(\bar U))-\Phi_0,\psi}
    &=\ip{\Phi_0(s)+w-\Phi_0(\cdot-\gamma_0(\bar U)),\psi(\cdot-\gamma_0(\bar U))}.
\end{align*}

We now assume that $\gamma_0$ has a first and second order Fr\'echet derivative, meaning that we can expand $\gamma_0$ around $\Phi_0(s)$ in the following way:
\begin{align}
    \gamma_0(\bar U)=\pi_*(s)+D\gamma_0(\Phi_0(s))[w]+\frac12D^2\gamma_0(\Phi_0(s))[w,w]+\mathcal{O}(w^3).
\end{align}
Upon plugging this back into the equation for $\gamma_0$, we find
\begin{align*}
    \ip{\Phi_0(s)&+w-\Phi_0(s)--\Phi'_0(s)D\gamma_0(\Phi_0(s))[w]+\frac12\left(\Phi'_0(s)D^2\gamma_0(\Phi_0(s))[w,w]-\Phi''_0(s)\left(D\gamma_0(\Phi_0(s))[w]\right)^2\right)\\
    &+\mathcal{O}(w^3),\psi(s)-\psi'(s)D\gamma_0(\Phi_0(s))[w]+\mathcal{O}(w^2)}=0.
\end{align*}
Collecting terms of equal order in $w$, we find at $\mathcal{O}(1)$ that $\ip{0,\psi(s)}=0$ which is automatically satisfied. For the next order, we find
\begin{align*}
    \ip{w+\Phi'_0(s)D\gamma_0(\Phi_0(s))[w],\psi(s)}=0.
\end{align*}
Hence
\begin{align*}
    D\gamma_0(\Phi_0(s))[w]=-\ip{w,\psi(s)}.
\end{align*}
When we collect terms of $\mathcal{O}(w^2)$, we find
\begin{align}
\begin{split}
    \frac12 D^2\gamma_0(\Phi_0(s))[w,w]&-\frac12\ip{\Phi_0''(s),\psi(s)}\left(D\gamma_0(\Phi_0(s))[w]\right)^2\\&-\ip{\Phi_0'(s),\psi'(s)}\left(D\gamma_0(\Phi_0(s))[w]\right)^2-\ip{w,\psi'(s)}D\gamma_0(\Phi_0(s))[w]=0.
    \end{split}
\end{align}
Plugging in the result for the first derivative and solving the equation, then results in
\begin{align*}
    D^2\gamma_0(\Phi_0(s))[w,w]=-\ip{\Phi''_0,\psi}\ip{w,\psi(s)}^2-2\ip{w,\psi'(s)}\ip{w,\psi(s)}.
\end{align*}

\subsection{Fr\'echet derivative of orbital drift}
We define 
\begin{align}
    J(\bar U)=\int_0^\infty a(v(s'))ds'.
\end{align}
Now, note that we do not take the derivative with respect to $v(s')$, but with respect to a direction $w$ which is related to the initial condition of $v(s')$. We explicitly write $\gamma_0(\Phi_0(s),w)$ to indicate that we introduce a perturbation $w$ away from $\Phi_0(s)$. Hence,
\begin{align}
\begin{split}
v(0):=v_0&=\bar U(\cdot+\gamma_0(\Phi_0(s),w))-\Phi_0\\
&=\Phi_0(\cdot-\pi(s)+\gamma_0(\Phi_0(s),w))+w(\cdot+\gamma_0(\Phi_0(s),w))-\Phi_0. 
\end{split}
\end{align}
Note that for $w=0$, $\gamma_0(\Phi_0(s),0)=\pi_*(s)$, hence $v_0=0$ as expected. Given the fact that 
\begin{align*}
    a(v(s'))=\frac{\ip{N(v(s')),\psi}}{1-\ip{v(s'),\psi'}}=\frac12\ip{f''(\Phi_0)v(s')^2,\psi}+\mathcal{O}(v(s'))^3,
\end{align*}
we find that the second-order Fr\'echet derivative of $a$ with respect to $v(s')$ is given by $\ip{f''(\Phi_0)\cdot^2,\psi}$. As the first Fr\'echet derivative $a(v(s')$ is zero, we therefore find that, by the chain rule,
\begin{align}
    D^2a(v(s'))[w,w]=D^2a(\Phi_0(s))[Dv(s')[w],Dv(s')[w]].
\end{align} 
Hence, we are only left with the Fr\'echet derivative of $v(s')$. Given the fact that $v(s')$ is the solution of a nonlinear PDE with initial condition $v(0)$, the first Fr\'echet derivative is the linearized equation with as initial condition the linearization of $v(0)$ around $v$. Hence:
\begin{align}
    Dv(s')(\Phi_0(s))[w]=S(s')[w(\cdot+\pi_*(s))-\Phi_0'\ip{w,\psi(\cdot-\pi_*(s))}].
\end{align}
Working our way back up again, we now find that 
\begin{align}
    D^2J(\Phi_0(s))[w,w]=\int_0^\infty\ip{f''(\Phi_0)\left(S(s')[w(\cdot+\pi_*(s))-\Phi_0'\ip{w,\psi(\cdot-\pi_*(s))}]\right)^2,\psi}ds'.
\end{align}

\subsection{Computing infinite sums}
In order to compute the autonomous isochronal phase, we need to plug the Fr\'echet derivatives back into the It\^o lemma ... Hence, we have to compute infinite sums. Of course, we could just leave these expressions, but we can wonder if we can find closed expressions. As a start, we study the lowest-order expansion of the phase, as we have a single expression for its variance. 
When we define 
\begin{align*}
    I(t)=\int_0^t\ip{g(\Phi_0)dW^Q_s(\cdot+\pi_*(s)),\psi}, 
\end{align*}
we can apply an It\^o lemma to find that
\begin{align}
 E[I^2(t)]=\sum_{i=0}^\infty\int_0^t \ip{g(\Phi_0)T_{-\pi_*(s)}\sqrt{Q}e_i,\psi}^2ds.
\end{align}
Now, by translation invariance, we can move the translation through the $\sqrt{Q}$ and find
\begin{align*}
 E[I^2(t)]&=\sum_{i=0}^\infty\int_0^t \ip{g(\Phi_0)\sqrt{Q}T_{-\pi_*(s)}e_i,\psi}^2ds.\\
 &=\int_0^t \nrm{\ip{g(\Phi_0)\sqrt{Q}T_{-\pi_*(s)}\cdot,\psi}}_{HS(L^2,\R)}^2ds.
\end{align*}
Now, note that for each value of $\pi_*(s)$, the set $(T_{-\pi_*(s)}e_i)_i$ is an orthonormal basis of $L^2(\R)$ when $(e_i)_i$ is an orthonormal basis. The Hilbert-Schmidt norm is independent of the chosen basis, so we can now drop the translation and the integrand becomes time-independent, i.e.  
\begin{align}
 E[I^2(t)]&=\sum_{i=0}^\infty\ip{g(\Phi_0)\sqrt{Q}e_i,\psi}^2t.
\end{align}
Next, we write out the $L^2$-innerproduct and the convolution $\sqrt{Q}e_i$, resulting in a total of four integrals:
\begin{align*}
    \sum_{i=0}^\infty\ip{g(\Phi_0)\sqrt{Q}e_i,\psi}^2&=\sum_{i=0}^\infty\int_\R\int_\R g(\Phi_0(x))\sqrt{Q}e_i(x)\psi(x)g(\Phi_0(x'))\sqrt{Q}e_i(x')\psi(x') dx'dx\\
    &=\sum_{i=0}^\infty\int_\R\int_\R g(\Phi_0(x))\psi(x)g(\Phi_0(x'))\psi(x')\int_\R p(x-y)e_i(y)dy\int_\R p(x'-y')e_i(y')dy' dx'dx.\\
\end{align*}
By Plancherel's theorem, we know that
\begin{align*}
   \sum_{i=0}^\infty \int_\R p(x-y)e_i(y)dy\int_\R p(x'-y')e_i(y')dy'
   &=\int_\R p(x-z)p(x'-z)dz\\
   &=\int_\R p(x-x'+y)p(y)dy\\
   &=\int_\R p(x-x'-y)p(y)dy:=q(x-x')
\end{align*}
where the switch $y\to -y$ followed from the fact that $p$ is symmetric. Plugging everything back together now results in
\begin{align*}
   \sum_{i=0}^\infty\ip{g(\Phi_0)\sqrt{Q}e_i,\psi}^2
    &=\int_\R\int_\R g(\Phi_0(x))\psi(x)g(\Phi_0(x'))\psi(x')q(x-x')dx'dx\\
    &=\ip{g(\Phi_0)\psi,Qg(\Phi_0)\psi},
\end{align*}
as predicted in the introduction. Also note that from this expression, it is immediately clear that the sum does not depend on any shift. Hence, the abstract statement that the sum must be independent of the shift because the HS-norm is independent of the chosen basis, can now be seen directly. 

With these computations under our belt, we can now compute the average autonomous isochronal phase. For the second-order derivative of $\gamma_0$ we now find
Plugging this back into the It\^o formula results in
\begin{align*}
    \frac12&\sum_{i=0}^\infty D^2\gamma(\Phi_0(s))[g(\Phi(s))T_{-\pi_s(s)}\sqrt{Q}e_i,g(\Phi(s))T_{-\pi_s(s)}\sqrt{Q}e_i]\\&=-\frac12\sum_{i=0}^\infty\ip{\Phi_0'',\psi}\ip{g(\Phi(s))T_{-\pi_s(s)}\sqrt{Q}e_i,\psi(s)}^2-\ip{g(\Phi(s))T_{-\pi_s(s)}\sqrt{Q}e_i,\psi'(s)}\ip{g(\Phi(s))T_{-\pi_s(s)}\sqrt{Q}e_i,\psi(s)}\\
    &=-\frac12\ip{\Phi_0'',\psi}\ip{Qg(\Phi_0)\psi,g(\Phi_0)\psi}-\ip{Qg(\Phi_0)\psi',g(\Phi_0)\psi}.
\end{align*}
The first sum immediately follows from the computation above, the next one follows after replacing one of the $\psi$ with $\psi'$.

\end{document}